\documentclass[portuges,12pt,letter]{article}
\usepackage[centertags]{amsmath}
\usepackage{amsfonts}
\usepackage{newlfont}
\usepackage{amscd}
\usepackage{graphics}
\usepackage{epsfig}
\usepackage{indentfirst}
\usepackage{amsxtra}
\usepackage[latin1]{inputenc}
\usepackage{amssymb, amsmath}
\usepackage{amsthm}
\usepackage[mathscr]{eucal}

\title{\bf Optimal control for the nonlinear Fisher-Kolmogorov system with applications to aquatic plant management}%\thanks{Grants or other notes

\author{Alexandre Molter$^1$ and Fabio Silva Botelho$^2$ \\ \\
$^1$Mathematics and Statistics Department \\
Federal University of Pelotas - UFPel \\
Pelotas, RS - Brazil \\ \\
$^2$ Mathematics Department \\
Federal University of Santa Catarina - UFSC \\
Florian\'{o}polis, SC - Brazil
}
\date{}
\begin{document}
\maketitle

\abstract{ Spatiotemporal dynamics of populations may be described by the reaction-diffusion Fisher-Kolmogorov model. In this work we have proposed a new formulation for a control problem of aquatic plants in a temporal dynamics. The solution of this problem is extended to a spatiotemporal Fisher-Kolmogorov system with multiple species of plants interacting in the same place. The control consists on human intervention as a strategy for management of the aquatic plants. In our applications, one plant and  two plants cases  have been considered. Simulation results are presented to show the effectiveness of the proposed control strategies.
 }

\section{Introduction}
\label{intro}

The population control for biological systems is currently used in several areas of knowledge, to control for example, the behavior of some animal, plant, bacteria and cell systems, or inhibit the growth of the concerning populations. Thus, a mathematical formulation for biological control problems may be helpful to describe the behavior of these populations and make appropriate predictions.

Mathematical formulations for biological control strategies have been presented in \cite{Raf,Rafi,Tan,Goh,Mol} where it is indicated  this kind of control, as applicable, may cause less  environment damages. Different forms of control have been so far used, such as an optimal nonlinear control  and respective feedback, optimal linear control for nonlinear systems and the use of Lyapunov functions. When dealing with aquatic plants, three forms of control are possible: the control through the use of chemicals,  a  physical control or a biological one (please see \cite{Osmond}, for instance). In the present work, we consider a physical control through human management.

The use of water hyacinths as an alternative to waste-water treatment is presented in \cite{Per,Kawai,Nest,Mishra,Saha}, and \cite{Hen,Henry}, showing interactions of two or more aquatic plants using  a waste-water treatment. The works \cite{Costa,Raf} solve optimization problems to find the optimal population level which keeps a best performance of the concerning system, taking into account only a temporal evolution. Spatiotemporal control has been proposed in the recent works \cite{Daiane}, used as an application to a water hyacinths population.

Concerning the nonlinear systems with reaction-diffusion, such as the Fisher-Kolmogorov \cite{Mur,Log}, there are several works which seek to show the coexistence between species and obtain concerning solutions  \cite{Bra,Halla,Vene,Ale,Muha,Vee,Kam}. In these works, theoretical studies on equilibrium and the existence of solutions are usually presented and relating numerical solutions are also developed. As a study on the control of Fisher-Kolmogorov systems we highlight the work of \cite{Duan}, however, in this work there is no applications, but only the formulation and proof of some theorems of the existence of optimal solutions for concerning control problems. Finding the control functions for the nonlinear Fisher-Kolmogorov systems may  be a complex task. Indeed there are only a few works on this type of control.

Based on \cite{Costa,Daiane,Raf} this work proposes a new formulation for the control problem of aquatic plants in stabilization ponds. Initially, a temporal optimization problem is solved, using a general Lotka-Volterra model. Afterwards, this formulation is extended to a spatiotemporal system with multiple species of plants interacting in the same pond, using a general Fisher-Kolmogorov model. The control consists on human intervention as a strategy for management of the aquatic plant population in polluted waters. As applications, in a first step we consider one plant in the pond and in a second step two plants in the pond are considered. Simulations are developed to show the efficiency of such  proposed control strategies.

This work is organized as follows. In the Section 2, the mathematical modeling of the control problem is presented. In the Section 3, we show the equation solution techniques. Further, in Section 4, applications to aquatic plants management are developed and the results discussed. Finally, some conclusions and remarks are presented.

%%%%%%%%%%%%%%%%%%%%%%%%%%%%%%%%%%%%%%%%%%%%%%%%%%%%%%%%%%%%%%%%%%%%%%%%%%%%%%%%%%%%%%%%%%%%%%%%%%%%%%%%%%%%%
\section{Mathematical modeling of the control problem}
\label{sec:1}
Consider a general Lotka-Volterra model of $N$ interacting populations described by
\begin{equation}\label{eq1}
\frac{d w_j}{dt}= a_j w_j-\sum_{k=1}^N b_{jk}w_jw_k, \;\; \forall j \in \{1,\cdots,N\},
\end{equation}
where $w_j(t)$ is the density of population $j$ at the instant $t$; $a_j$ is the growth rate of population $j$ , $b_{jk}$ are the competition
coefficients between the populations.

We consider the action of a control variable which suppress a certain quantity of population from the system.  Such a control variable is point-wise denoted by $u_j(t)$ and introduced in the equation system as it follows:
\begin{eqnarray}\label{eq2}
\frac{d w_j(t)}{dt}= a_j (w_j(t)-u_j(t))\nonumber \\
-\sum_{k=1}^N b_{jk}(w_j(t)-u_j(t))(w_k(t)-u_k(t))\nonumber \\
-\tau_j u_j,\; \forall j \in \{1,\cdots,N\},
\end{eqnarray}
where $\tau_j$ are positive constants that characterize the technical conditions of population withdraw. Equation (\ref{eq2}) describes the dynamics
of the system (\ref{eq1}) with an application of the control $u_j(t)$.

Therefore, denoting by $[0,T]$ the time interval of the control action, respectively, we have
\begin{eqnarray}w_j(T)=w_j(0)+\int_0^T \frac{d w_j(t)}{dt}\;dt \nonumber \\
= w_j(0)+\int_0^T\left( a_j (w_j(t)-u_j(t))\right. \nonumber \\
 \left.-\sum_{k=1}^N b_{jk}(w_j(t)-u_j(t))(w_k(t)-u_k(t))-\tau_j u_j\right)\;dt.
\end{eqnarray}

Consider now the problem of minimizing
\begin{equation}J(w,u)= -\sum_{j=1}^N w_j(T)-\int_0^T \tau_j u_j(t) \;dt, \end{equation}
with the constraints
\begin{equation}0\leq u_j(t) \leq w_j(t),\; \text{on}\;\; [0,T],\; \forall j \in \{1,\cdots,N\}.\end{equation}

Thus, we shall look for a critical point of the Lagrangian

\begin{eqnarray}&& L(w,u,\lambda)=J(w,u) +\nonumber \\ &&
 \sum_{k=1}^N\left(\int_0^T \lambda_1^k(t)(u_k(t)-w_k(t))\;dt+\int_0^T \lambda_2^k(t) u_k(t) \;dt\right)\nonumber \\
 &=& - \sum_{j=1}^Na_j \int_0^T (w_j(t)-u_j(t))\;dt + \nonumber \\ &&
 \sum_{j,k=1}^N b_{jk}\int_0^T(w_j(t)-u_j(t))(w_k(t)-u_k(t))\;dt + \nonumber \\ &&
 \sum_{j=1}^N\left(\int_0^T \lambda_1^j(t)(u_j(t)-w_j(t))\;dt \right. \nonumber \\ && \left.+\int_0^T \lambda_2^j(t)u_j(t)\;dt - w_j(0)\right).
\end{eqnarray}

The variation of $L$ in $u_k$, gives the stationary equation (see \cite{Botelho} for related theoretical optimization results)
\begin{equation}a_j +\sum_{k=1}^N (b_{jk}+b_{kj})(u_k-w_k)+\lambda_1^j+\lambda_2^j=0,\end{equation}
so that
\begin{equation}\{u_j-w_j\}=-\left\{ b_{jk}+b_{kj}\right\}^{-1} \left\{ a^k +   \lambda_1^k+\lambda_2^k\right\}.\end{equation}

For each $k \in \{1,\cdots,N\}$, we have two possibilities, either $\lambda_1^k(t)=0, \; \lambda_2^k(t)=0$ (in such a case the constraints are not active); or $\lambda_1^k(t)=0$ and $\lambda_2^k(t) >0,$ which corresponds to the constraint $u_k(t)=0$ to be active, so that in such a case we have a local minimum. The possibility
with $\lambda_1^k(t)>0$ and $\lambda_2^k(t) =0$ at any point $t$ means $u_k(t)=w_k(t)$ does not correspond to a local minimum.

Hence defining
\begin{equation}\label{xi}
\{\xi_k\}=\left\{b_{jk}+b_{kj}\right\}^{-1} \left\{ a_j\right\},
\end{equation}
and \begin{equation}\label{ak}
A_j(t)=\left(-a_j +\sum_{k=1}^N (b_{jk}+b_{kj}) w_k\right),
\end{equation}
 we obtain
\begin{equation}\label{cont}
 u_k(t)=\left\{
\begin{array}{ll}
 w_k(t)-\xi_k,& \text{ if } A_k(t)>0
 \\
 0,& \text{ if } A_k(t) \leq 0,\end{array} \right.
 \end{equation}
 $\forall k \in \{1,\cdots,N\}.$

The parameter $\xi_k$ in (\ref{xi}) corresponds to the optimal level of population, which means  the
control is applied when the amount of the population surpass this level.

Replacing the control function (\ref{cont}) into the equation (\ref{eq2}) yields
 $$\frac{d w_j(t)}{dt}=$$
\begin{equation}\label{siscont}
 \left\{
\begin{array}{ll}
\displaystyle a_j\xi_j -\sum_{k=1}^N b_{jk}(\xi_j)(\xi_k)-\tau_j(w_j(t)-\xi_j),\; & \text{ if } A_j(t)>0 \\
\displaystyle a_j w_j-\sum_{k=1}^N b_{jk}w_jw_k, \;\; & \text{ if } A_j(t) \leq 0,
 \end{array} \right.
 \end{equation}
 $\forall j \in \{1,\cdots,N\}.$

From (\ref{siscont}) we may observe that when the control is
applied, the differential equations system becomes linear due to some algebraic simplifications. We can also observe that the equation system (\ref{siscont}) have an exact analytic solution in each instant $t$, taken into account an initial condition $w_j(0)=(w_{0})_j$. In the works of \cite{Costa} and \cite{Daiane} it may be seen the application of such a system to the population of some plant species.

Next, consider the general model of Fisher-Kolmolgorov describing temporal and spatial dynamics of $N$ interacting populations
\begin{equation}\label{eqkol}
\frac{\partial w_j}{\partial t}= D_j\frac{\partial^{2} w_j}{\partial x^{2}}+a_jw_j-\sum_{k=1}^N b_{jk}w_jw_k,
\end{equation}
$\forall j \in \{1,\cdots,N\},$\\
where $D_j$ is the coefficient of diffusion for each species and $w_j=w_j(x,t)$.

Let $u_k(x,t)$ be the control function to be introduced in the equation system (\ref{eqkol}), which satisfies the same conditions as $u_k(t)$ previously defined, but now considering also a dispersal in the space  and $\xi_k$ with same dimension as $w_k$. The new control function can be taken in the form
\begin{equation}\label{contspace}
 u_k(x,t)=\left\{
\begin{array}{ll}
 w_k(x,t)-\xi_k,& \text{ if } A_k(x,t)>0
 \\
 0,& \text{ if } A_k(x,t) \leq 0,\end{array} \right.
 \end{equation}
 $\forall k \in \{1,\cdots,N\}.$

 Replacing the control function (\ref{contspace}) in equation (\ref{eqkol}), similarly as in (\ref{eq2}), yields:
\begin{equation}\label{edpcont}
\frac{\partial w_j}{\partial t}= \left\{
\begin{array}{cc}
\displaystyle D_j \frac{\partial^{2} w_j}{\partial x^{2}}+a_j (w_j-u_j)\\
\displaystyle -\sum_{k=1}^N b_{jk}(w_j-u_j)(w_k-u_k)\\
\displaystyle -\tau_j u_j,\; \quad \quad \quad \text{ if } A_j(x,t)>0,\\
\begin{small}\displaystyle D_j \frac{\partial^{2} w_j}{\partial x^{2}}+a_jw_j-\sum_{k=1}^N b_{jk}w_jw_k,\text{ if } A_j(x,t) \leq 0.\end{small}
\end{array}
\right.
\end{equation}
$\forall j \in \{1,\cdots,N\}.$

Through appropriate algebraic manipulations in the first equation of the system (\ref{edpcont}), we may obtain
\begin{equation}\label{edpcont2}
\frac{\partial w_j}{\partial t}= \left\{
\begin{array}{cc}
\displaystyle D_j \frac{\partial^{2} w_j}{\partial x^{2}}+a_j (\xi_j)-\sum_{k=1}^N b_{jk}(\xi_j)(\xi_k)\\
\displaystyle -\tau_j(w_k-\xi_k),\quad \quad \quad \text{ if } A_j(x,t)>0,\\
\begin{small}\displaystyle D_j \frac{\partial^{2} w_j}{\partial x^{2}}+a_jw_j-\sum_{k=1}^N b_{jk}w_jw_k,\text{ if } A_j(x,t) \leq 0.\end{small}
\end{array}
\right.
\end{equation}
$\forall j \in \{1,\cdots,N\}.$

 The initial and boundary conditions  depend on  the specificity of each application. For example, in an application of an aquatic plant system \cite{Daiane} the initial conditions may be given by:
\begin{equation}\label{condiinicial}
\displaystyle w_j(x,0)=(w_{0})_j,
\end{equation}
and the boundary conditions are of Neumann's type, that is
\begin{equation}\label{condicoesw}
\begin{array}{ccc}%\label{sistemacomcontrole}
\displaystyle \frac{\partial w_j(0,t)}{\partial x}=0; &  \displaystyle\frac{\partial w_j(l,t)}{\partial x}=0, &  \ \ t \geq 0,
\end{array}
\end{equation}
which indicates  there is no flow of plants at the pond borders (aquatic plants depend on water to survive).

Note that the first equation in (\ref{edpcont2}) is similar as the first equation in (\ref{siscont}), considering that $w(t)$ is replaced by $w(x,t)$  and that it has been added $D_j$, and a diffusive term $w_{xx}$. Moreover, the first equation in (\ref{edpcont2})  becomes linear with the addition of the control strategy. In the next section we show a technique to solve equation systems such as indicated in (\ref{edpcont2}), first and second equations.

%%%%%%%%%%%%%%%%%%%%%%%%%%%%%%%%%%%%%%%%%%%%%%%%%%%%%%%%%%%%%%%%%%%%%%%%%%%%%%%%%%%%%%%%%%%%%%%%%%%%%
\section{Equations solution technique }

%%\subsection{The algorithm}
In this section we present the numerical algorithm, in a finite differences and numerical Newton's Method similar context (not exactly the Newton's method but an appropriate adaptation), to be used in the subsequent section for the solution of the partial differential equations given in the mathematical formulation of each problem considered.

At this point we denote the matrix $D_{xx}$ by

\begin{gather}\label{910} D_{xx}=\left[
\begin{array}{cccccc}
-1 & 1&0 & 0 & \cdots & 0 \\
1 & -2 & 1&0 & \cdots & 0 \\
0 & 1 &-2 & 1 & \cdots & 0 \\
\vdots & \vdots &    \vdots &\ddots & \vdots & \vdots \\
0 & 0 & \cdots & 1 & -2 & 1 \nonumber \\
0& 0 &  \cdots & 0 & 1 & -1
\end{array} \right]/d_1^2,\end{gather}
where $d_1$ will be specified in the next lines.

Such an algorithm is to obtain a numerical solution on the set $[0,10] \times [0,T],$ where $[0,T]$ is the time interval considered.
\begin{enumerate}

\item Set $T=30 (days)$,  $N_1=5000$ (number of nodes in $[0,T]$), $N_2=100$ (number of nodes in $[0,10]$), $d=T/N_1$ $d_1=1/N_2$;

\item Set $\{w^0_j\}=w_0^j$;

\item For  $l=1:N_1$ do

\begin{itemize}
\item $p=1$,
\item $\{w_j^p\}=\{w_j^{l-1}\},\; \forall j \in \{1, \ldots,N\};$
\item $b_{12}=1.0$;
\item While $(b_{12}>10^{-4})$ and  ($p<100$) do
\begin{itemize}
\item For $j=1:N$ do
\begin{itemize}
\item For $s=1:N_2$ do

  $(A_j)_s=(-a_j+\sum_{k=1}^N(b_{jk}+b_{kj})(w_j^p)_s);$

  $(u_j^s)=0;$

  If $(A_j)_s >0$, then $(u_j)_s=(w_j^p)_s-\xi_j$;
  \item end;
  \end{itemize}
  \item end;
  \end{itemize}

  \item Calculate $\hat{w}$ by solving the linearized equation

  $$\frac{ \hat{w}_j-w_j^{l-1}}{d}= D_j D_{xx} \hat{w}_j+H_j(\hat{w}),\; \forall j \in \{1, \ldots, N\};$$
  where,
  \\
  If $(A_j)_s\leq0$, then

   $(H_j)_s(\hat{w})= a_j (\hat{w}_j)_s-2\sum_{k=1}^N b_{jk} (\hat{w}_j)_s (w_k^p)_s+\sum_{k=1}^N b_{jk}(w_j^p)_s(w_k^p)_s;$
\\
\\
  If $(A_j)_s> 0$, then

  $(H_j)_s(\hat{w})= a_j \xi_j-\sum_{k=1}^N b_{jk}\xi_j \xi_k -\tau_j (u_j)_s;$
  \\
  \item Set $B_j=\|\hat{w}_j-w_j^p\|_\infty, \; \forall j \in \{1,\ldots,N\};$
  \item  Set $b_{12}=\max_{ j \in \{1,\ldots,N\}} B_j;$
  \item $w_j^{p+1}=\hat{w}_j,\; \forall j \in \{1,\ldots,N\};$
  \item $p =p+1;$
  \item end;
    \item $w_j^l=w_j^p, \; \forall j \in \{1, \ldots, N\};$
    \end{itemize}
  \item End.
    \end{enumerate}
%\begin{remark} For the numerical results purposes we have  considered a simplified  control acting  on all $x \in [0,10]$ or not acting at any $x \in [0,10]$ at each concerning time $t$, which is slightly different from the above algorithm. \end{remark}
%%%%%%%%%%%%%%%%%%%%%%%%%%%%%%%%%%%%%%%%%%%%%%%%%%%%%%%%%%%%%%%%%%%%%%%%%%%%%%%%%%%%%%%%%%

\section{Applications to aquatic plants management}

For the applications we have considered a system with two aquatic plants, {\it Eichhornia crassipes} and {\it Pistia stratiotes}, often used in waste-water treatment. There is a threshold, namely $\xi$, to consider the growth of the two aquatic plants, where a greater accumulation of biomass corresponds to a greater capacity of extraction of residues, provided that the threshold is not exceeded, which will result in a loss of capacity to remove toxic products from the water.
Therefore, to maintain the pond efficiency, one have to remove plants daily, in order to hold high rates of constant growth, without
causing undesirable side effects. The amount of plants to be removed daily was established in the previous section, equation (\ref{cont}) for a temporal system and (\ref{contspace}) for a spatiotemporal system, considering the values of the coefficients given in the next lines.
For the applications considered, the space is a channel fulfilled with residual water.

The function $w(x,t)$ in this application represents the plants density (dry mass $gm^{-2}t^{-1}$) in the space $x$ ($m^2$) and at
the instant $t$ (days).

\subsection{Application using one plant}

 The interaction between a specie with its environment, may be described by the partial differential equation of Fisher-Kolmogorov \cite{Log,Mur}:
\begin{equation}\label{equafisher}
\frac{\partial w}{\partial t}= D\frac{\partial^{2} w}{\partial x^{2}}+aw-bw^2,
\end{equation}

Replacing the system (\ref{edpcont2}), with control,  in equation (\ref{equafisher}) yields
\begin{small}
\begin{equation}\label{edpcont3}
\frac{\partial w}{\partial t}= \left\{
\begin{array}{cc}
\displaystyle D\frac{\partial^{2}w}{\partial x^{2}}+a(w-u)-b(w-u)(w-u) - \tau u, \; \text{ if } A>0,\\
\displaystyle D\frac{\partial^{2}w}{\partial x^{2}}+aw-bw^2,\;\text{ if } A\leq0,
\end{array}
\right.
\end{equation}
\end{small}

Replacing $u=w-\xi$ into (\ref{edpcont3}) we obtain
\begin{equation}\label{edpcont4}
\frac{\partial w}{\partial t}= \left\{
\begin{array}{cc}
\displaystyle D\frac{\partial^{2}w}{\partial x^{2}}+a\xi-b\xi^2 -\tau w +\tau\xi, \; \text{ if } A>0,\\
\displaystyle D\frac{\partial^{2}w}{\partial x^{2}}+aw-bw^2,\;\text{ if } A\leq0,
\end{array}
\right.
\end{equation}

Note that if we take $A=0$ and $b_{jk}$=$b_{kj}$, with a single type of plant, from equations (\ref{xi}) and (\ref{ak}) we may calculate
\begin{equation}
-a +(2b)w=0 \;\Rightarrow w(x,t)=\frac{a}{2b}.
\end{equation}
So, $A>0 \Rightarrow w(x,t)>\dfrac{a}{2b}$ and $A\leq0 \Rightarrow w(x,t)\leq\dfrac{a}{2b}$, which is consistent with the results of works \cite{Costa} and \cite{Daiane}.

Equation (\ref{edpcont4}) can be rewritten as
\begin{equation}\label{edpcont5}
\frac{\partial w}{\partial t}= \left\{
\begin{array}{cc}
\displaystyle D\frac{\partial^{2}w}{\partial x^{2}}-\tau w +P, \; \text{ if } A>0,\\
\\
\displaystyle D\frac{\partial^{2}w}{\partial x^{2}}+aw-bw^2,\;\text{ if } A\leq0,
\end{array}
\right.
\end{equation}
where $P = a^{2}/(4b)+a\tau/(2b)$ is the amount of plants in the channel with residual considering the control application, with is applied only when the plants population surpass the value $\xi=a/2b$. On the other hand, from equation (\ref{contspace}) we may obtain the amount of plants to be removed daily, $u^{\ast}(x,t)$, to keep the population at level $P$, being $u^{\ast}(x,t)=a^2/4b$.

The initial and boundary conditions have been taken as in (\ref{condiinicial}) and (\ref{condicoesw}), respectively.

\begin{figure}[!htb]
\begin{center}
	\includegraphics[height=10cm]{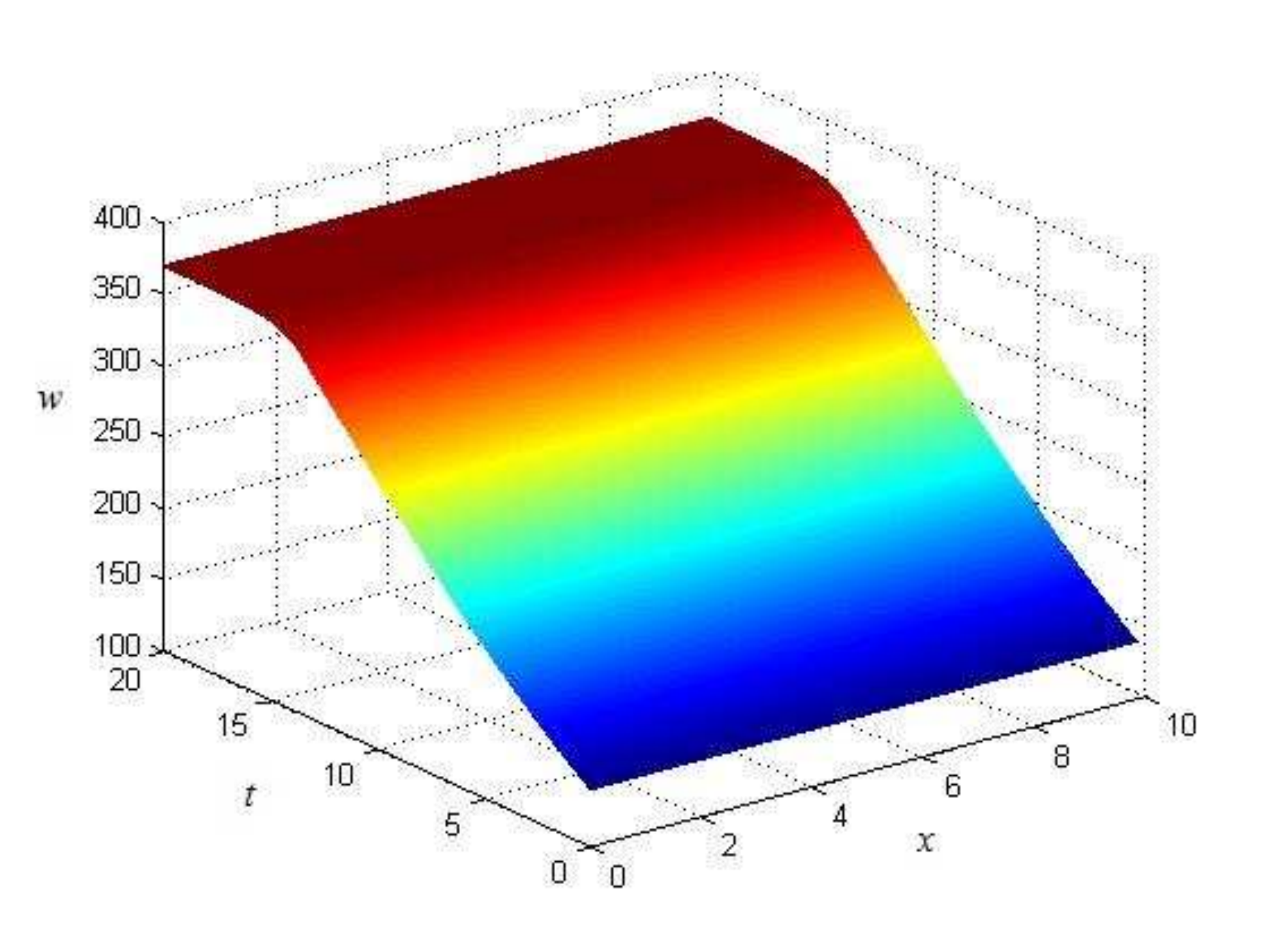}
    \vspace{-0.5cm}
\caption{\small{Solution $w(x,t)$ for the initial conditions $w(x,0)= 140$.}} \label{figure-Kol-5}
\end{center}
\end{figure}

\begin{figure}[!htb]
\begin{center}
	\includegraphics[height=10cm]{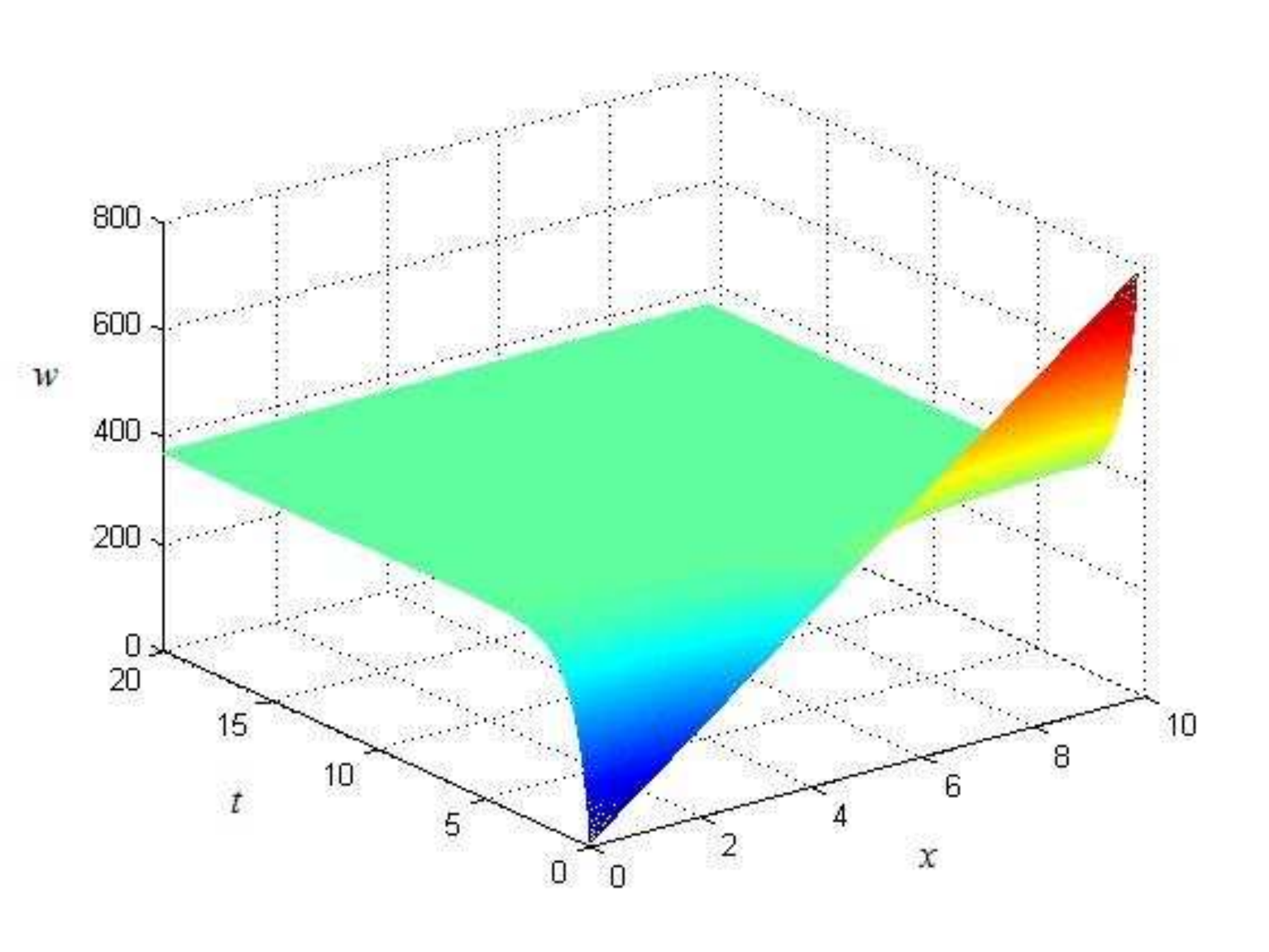}
    \vspace{-0.5cm}
\caption{\small{Solution $w(x,t)$ for the initial conditions $w(x,0)=80x$.}} \label{figure-Kol-6}
\end{center}
\end{figure}

The plant considered in this application is the {\it Eichhornia crassipes} and the related coefficients values have been taken from the works \cite{Costa,Daiane,Kawai,Per}. The coefficients values are: $K=700.68$, $a=0.103$, $b=0.000147$, $D=1.33$ and $\tau=1$. For the initial conditions
$w(x,0)=140$ and $w(x,0)=80x$, where $ x \in [0,10]$, the figures \ref{figure-Kol-5} and \ref{figure-Kol-6} show the behavior of $ w(x,t)$ considering a 20-day time evolution.

%\vspace{-6cm}
 We may observe in the figures \ref{figure-Kol-5} and \ref{figure-Kol-6} that in a few days, the population $w(x,t)$ is driven to the value of $P$, regardless if the number of plants is below or above this value.

\subsection{Application using two plants}

The of Fisher-Kolmogorov equations for two species is given by \cite{Bra,Mur}:
\begin{equation}\label{fisherkol}
\left\{
\begin{array}{c}
\displaystyle\frac{\partial w_1}{\partial t}= D_1\frac{\partial^{2} w_1}{\partial x^{2}}+a_1w_1-b_{11}w^2_1-b_{12}w_1w_2,\\
\\
\displaystyle\frac{\partial w_2}{\partial t}= D_2\frac{\partial^{2} w_2}{\partial x^{2}}+a_2w_2-b_{21}w_2w_1-b_{22}w^2_2.
\end{array} \right.
\end{equation}

Including a control in the system (\ref{fisherkol}), as shown in system (\ref{edpcont2}), yields
$$\frac{\partial w_1}{\partial t}=$$
\begin{small}
\begin{equation}\label{sistwo1}
\left\{
\begin{array}{c}
\displaystyle D_1 \frac{\partial^{2} w_1}{\partial x^{2}}+a_1 \xi_1-\sum_{k=1}^2 b_{1k}\xi_1\xi_k
 -\tau_1w_1+\tau_1\xi_1, \text{ if } A_1>0,\\
\displaystyle D_1\frac{\partial^{2} w_1}{\partial x^{2}}+a_1w_1-b_{11}w^2_1-b_{12}w_1w_2,\; \text{ if } A_1\leq0,\\
\end{array} \right.
\end{equation}
\end{small}
$$\frac{\partial w_2}{\partial t}=$$
\begin{small}
\begin{equation}\label{sistwo2}
\left\{
\begin{array}{c}
\displaystyle D_2 \frac{\partial^{2} w_2}{\partial x^{2}}+a_2 \xi_2-\sum_{k=1}^2 b_{2k}\xi_2\xi_k
 -\tau_2w_2+\tau_2\xi_2, \text{ if } A_2>0,\\
\displaystyle D_2\frac{\partial^{2} w_2}{\partial x^{2}}+a_2w_2-b_{21}w_2w_1-b_{22}w^2_2,\; \text{ if } A_2\leq0,
\end{array} \right.
\end{equation}
\end{small}

From equation (\ref{xi}) the values of $\xi_1$ and $\xi_2$ can be computed as:
\begin{equation}\label{xi2}
\left\{ \begin{array}{c}
\xi_1 \\
\xi_2 \\
\end{array} \right\}
=\left\{
\begin{array}{cc}
2b_{11} & b_{12}+b_{21} \\
b_{21}+b_{12} & 2b_{22} \\
\end{array} \right\}^{-1}
\left\{ \begin{array}{c}
a_1 \\
a_2 \\
\end{array} \right\},
\end{equation}
and equations (\ref{sistwo1}) and (\ref{sistwo2}) can be rewritten as
\begin{small}
\begin{equation}\label{sistwo3}
\frac{\partial w_1}{\partial t}=\left\{
\begin{array}{c}
\displaystyle D_1 \frac{\partial^{2} w_1}{\partial x^{2}}-\tau_1w_1+P_1, \text{ if } A_1>0,\\
\displaystyle D_1\frac{\partial^{2} w_1}{\partial x^{2}}+a_1w_1-b_{11}w^2_1-b_{12}w_1w_2,\; \text{ if } A_1\leq0,\\
\end{array} \right.
\end{equation}
\end{small}
\begin{small}
\begin{equation}\label{sistwo4}
\frac{\partial w_2}{\partial t}=\left\{
\begin{array}{c}
\displaystyle D_2 \frac{\partial^{2} w_2}{\partial x^{2}}-\tau_2w_2+P_2, \text{ if } A_2>0,\\
\displaystyle D_2\frac{\partial^{2} w_2}{\partial x^{2}}+a_2w_2-b_{21}w_2w_1-b_{22}w^2_2,\; \text{ if } A_2\leq0,
\end{array} \right.
\end{equation}
\end{small}
where
\begin{equation}\label{P1P2}
\left\{
\begin{array}{c}
\displaystyle P_1=a_1 \xi_1-b_{11}\xi_1\xi_1-b_{12}\xi_1\xi_2+\tau_1\xi_1 \\
\\
\displaystyle P_2=a_2 \xi_2-b_{21}\xi_2\xi_1-b_{22}\xi_2\xi_2+\tau_2\xi_2
\end{array} \right.
\end{equation}
are the amount of plants in the channel with residual water considering the control application, plant 1 and plant 2, respectively.

The initial and boundary conditions were taken in same form, for each specie, as in (\ref{condiinicial}) and (\ref{condicoesw}), respectively.
%%The initial conditions used are: $(w_1(x,0),w_2(x,0))=(280,80)$, $(w_1(x,0),w_2(x,0))=(700,350)$.
The technical capacity to remove the plants are given by $\tau_1=1$ and $\tau_2=1$ \cite{Costa,Daiane}. The coefficients of the system are those presented in the next lines, highlighting where they have been taken from or how they have been computed.

\noindent $\bullet$ {\it Plant 1 - Eichhornia crassipes}

\noindent $a_1=0.061$ \cite{Hen}; $K_1=993.2$ (in monoculture \cite{Hen});\\
$b_{11}= 0.0000614$ (computed from $K_1=a_1/b_{11}$);\\
$b_{12}=0.00001$ (estimate - according to \cite{Hen}, the interaction of the two species practically did not affect the water hyacinth growth);\\
$D_1=1.33$ (computed from the information containing in \cite{Per}).

\noindent $\bullet$ {\it Plant 2 - Pistia stratiotes}

\noindent $a_2=0.087$ \cite{Hen}; $K_2=436.7$ (in monoculture \cite{Hen});\\
$b_{21}=0.0001$ (estimate - according to \cite{Hen}, the interaction of the two species reduced significantly the Pistia growth);\\
$b_{22}=0.0001992$ (computed from $K_2=a_2/b_{22}$);\\
$D_2= 1.3$ (estimated to be the similar as water hyacinth, because its growth is similar).

\begin{figure}[!htb]
\begin{center}
	\includegraphics[height=10cm]{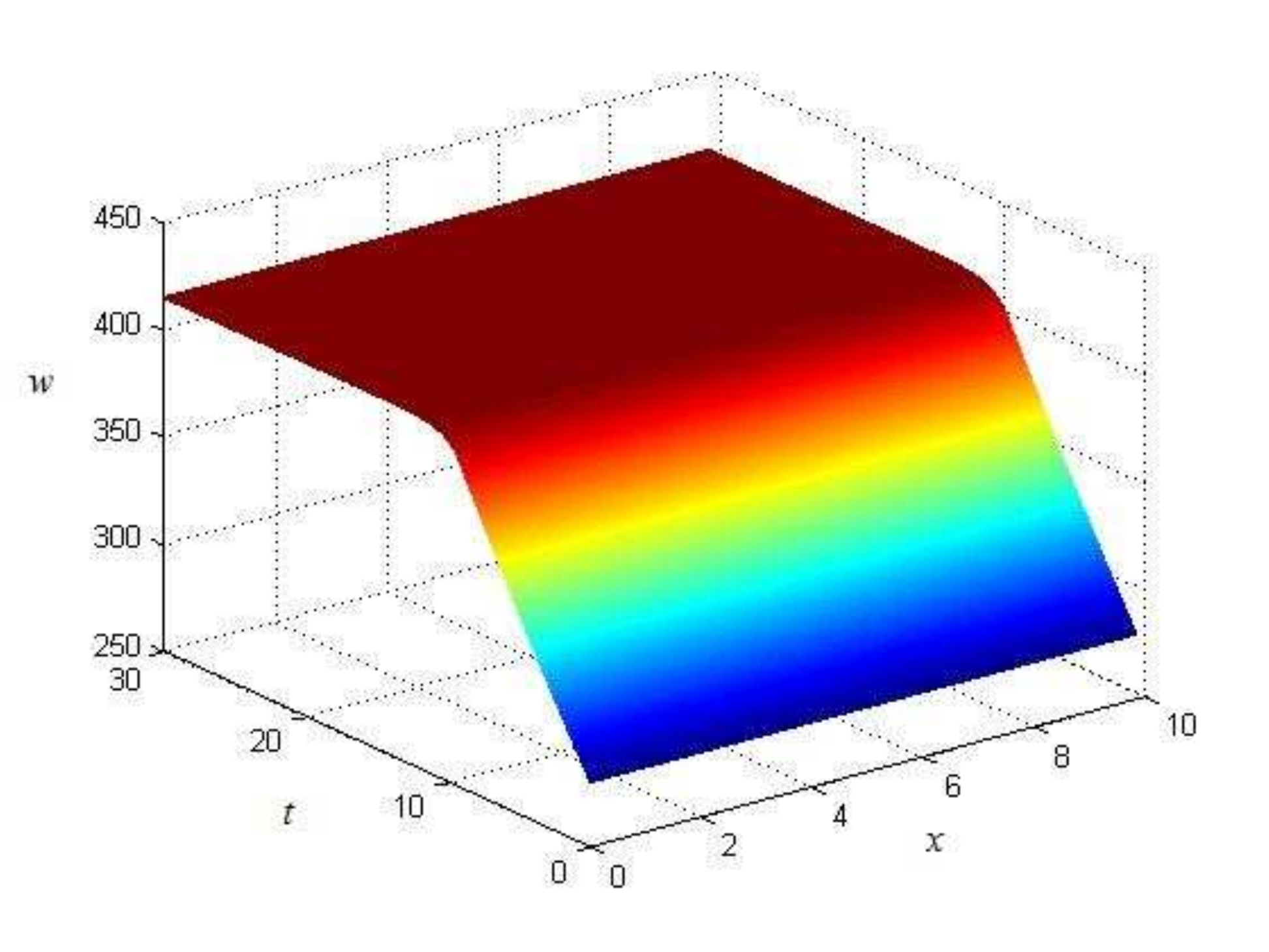}
    \vspace{-0.5cm}
\caption{\small{Solution $w_1(x,t)$ for the initial conditions $(w_1(x,0),w_2(x,0))=(280,80)$, plant 1.}} \label{figure-Kol-1}
\end{center}
\end{figure}

\begin{figure}[!htb]
\begin{center}
	\includegraphics[height=10cm]{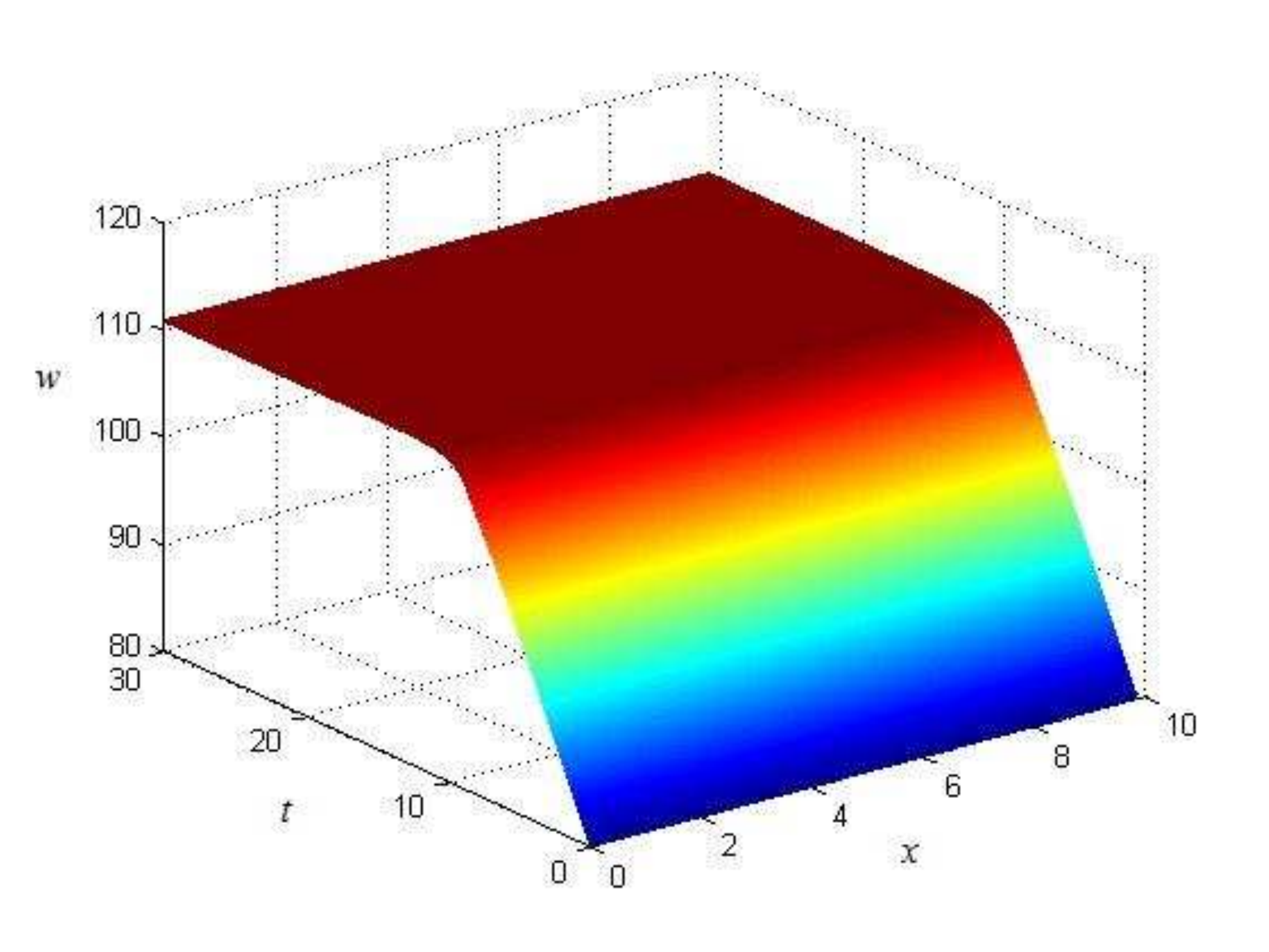}
    \vspace{-0.5cm}
\caption{\small{Solution $w_2(x,t)$ for the initial conditions $(w_1(x,0),w_2(x,0))=(280,80)$, plant 2.}} \label{figure-Kol-2}
\end{center}
\end{figure}

\begin{figure}[!htb]
\begin{center}
	\includegraphics[height=10cm]{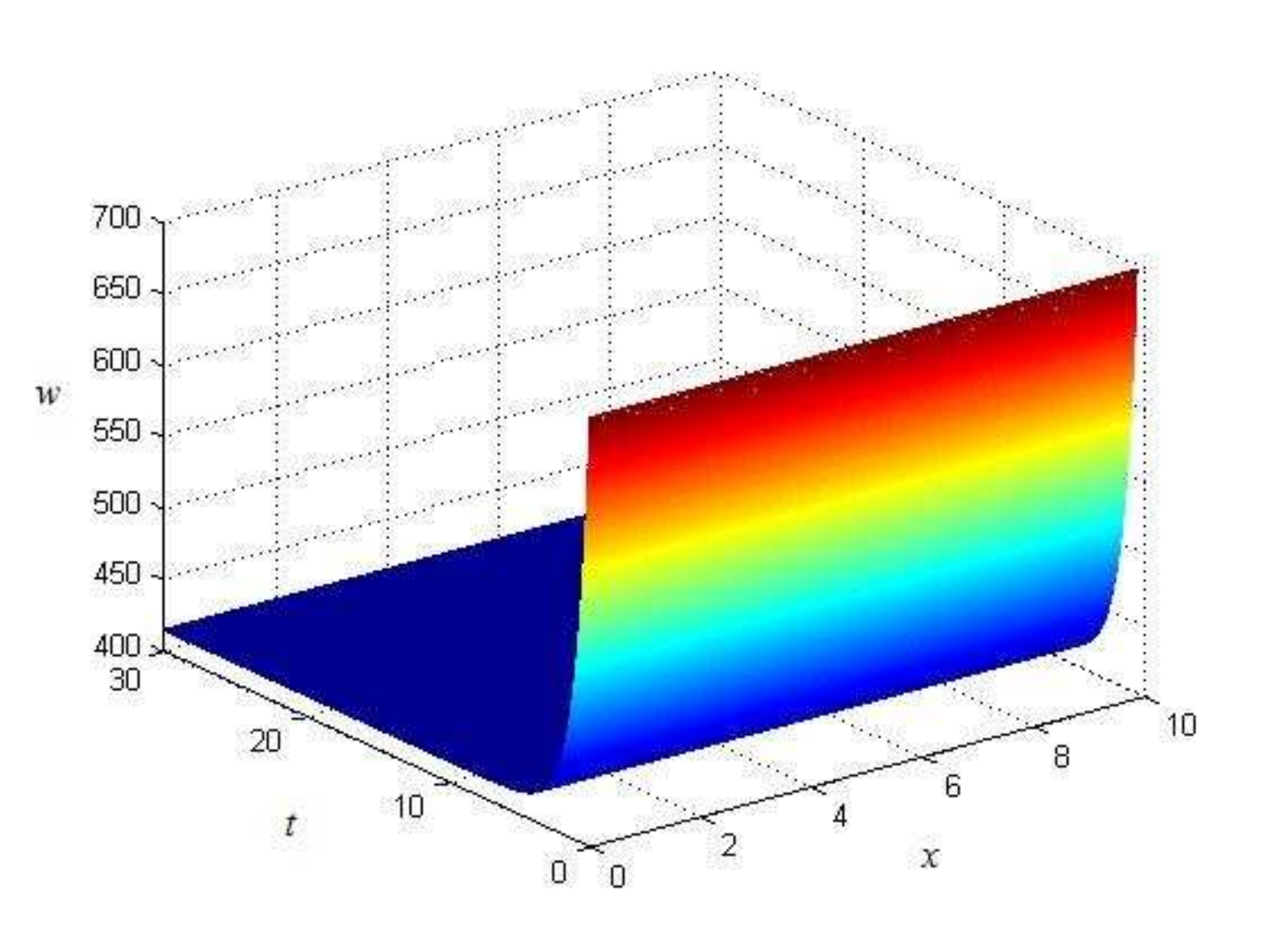}
    \vspace{-0.5cm}
\caption{\small{Solution $w_1(x,t)$ for the initial conditions $(w_1(x,0),w_2(x,0))=(700,350)$, plant 1.}} \label{figure-Kol-3}
\end{center}
\end{figure}

\begin{figure}[!htb]
\begin{center}
	\includegraphics[height=10cm]{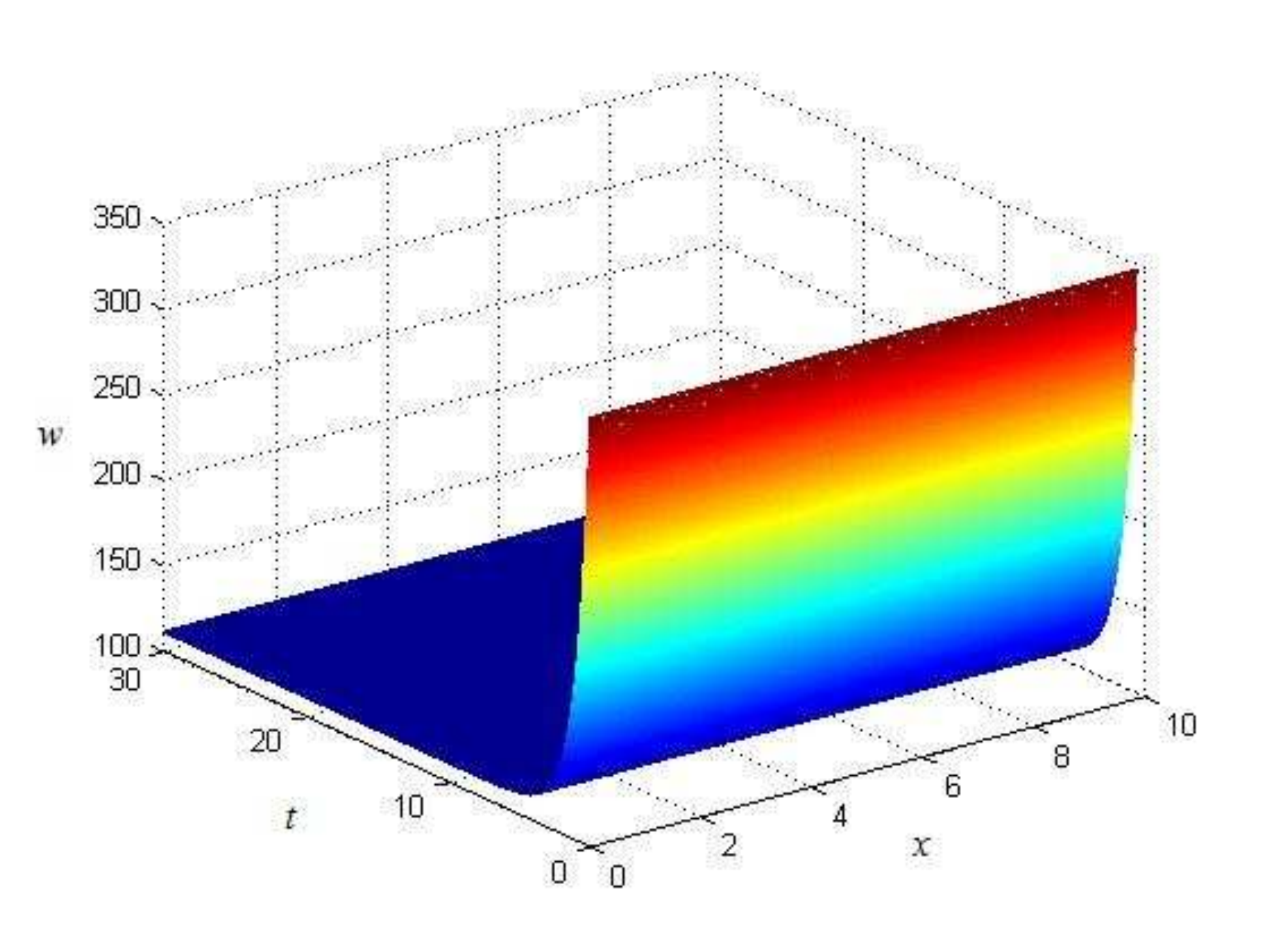}
    \vspace{-0.5cm}
\caption{\small{Solution $w_2(x,t)$ for the initial conditions $(w_1(x,0),w_2(x,0))=(700,350)$, plant 2.}} \label{figure-Kol-4}
\end{center}
\end{figure}

Figures \ref{figure-Kol-1} and \ref{figure-Kol-2} show the controlled growth of plants 1 and 2 into a $10m$ channel, i. e., $ x \in [0,10]$, and initial conditions $(w_1(x,0),w_2(x,0))=((280,80)$, respectively.
Figures \ref{figure-Kol-3} and \ref{figure-Kol-4} show the controlled growth of plants 1 and 2 into a $10m$ channel, i. e., $ x \in [0,10]$, and initial conditions $(w_1(x,0),w_2(x,0))=((700,350)$, respectively.
One can observe from figures \ref{figure-Kol-1}, \ref{figure-Kol-2}, \ref{figure-Kol-3}, and \ref{figure-Kol-4},  the system reaches the equilibrium level about 15 days when the initial condition is below the optimum level, and about 5 days when the initial condition is above the optimum level. The optimal level of controlled population is calculated from equation (\ref{P1P2}) and given by: $(w^{\ast}_1, w^{\ast}_2)=(414.23, 110.66)$. Thus, from equation (\ref{contspace}) can be calculated the amount of plants to be removed daily to maintain the desired level of plant population in such a system, namely, $(u^{\ast}_1, u^{\ast}_2)=(14.14, 2.75)$.

%%%%%%%%%%%%%%%%%%%%%%%%%%%%%%%%%%%%%%%%%%%%%%%%%%%%%%%%%%%%%%%%%%%%%%%%%%%%%%%%%%%

\section{Conclusions and remarks}

In this work a control has been applied to a system described by the Fisher-Kolmogorv model. An optimization problem was proposed and solved. The applications considered a system comprised by an aquatic plants population growing in a waste-water stabilization pond. Simulation results have been presented in which the proposed control has been efficient in maintaining the system at the optimal level of plants. In addition, from the applications presented, it may be inferred that:

\noindent $\bullet$ When populations of plants are below the optimal level, as time goes on, plants grow freely, in all directions, until they reach a  level for which the application of control begins, with daily withdrawal of plants.

\noindent $\bullet$ When populations of plants are intended above the optimal level, the control is immediately activated.

\noindent $\bullet$ To maintain the optimal level of {\it Eichhornia crassipes} in a channel with residual water $P=w^{\ast}=368.38$ one should remove $u^{\ast}=18.04$ daily (dry mass of plants).

\noindent $\bullet$ To maintain the optimal level of {\it Eichhornia crassipes} and {\it Pistia stratiotes} in a channel with residual water $(P_1,P_2)=(w^{\ast}_1, w^{\ast}_2)=(414.23, 110.66)$ one should daily remove $(u^{\ast}_1, u^{\ast}_2)=(14.14, 2.75)$ of plants (dry mass of each kind of plant).

The control considered in this work takes into account the removal of the plants being spread in the whole
pond area, ideally in each square meter. This  indicates us that stabilization ponds may be appropriately planned for a spatial
management.

\section*{Conflict of Interest}

The authors declare that they have no conflict of interest.

%\begin{acknowledgements}
%If you'd like to thank anyone, place your comments here
%and remove the percent signs.
%\end{acknowledgements}

% BibTeX users please use one of
%\bibliographystyle{spbasic}      % basic style, author-year citations
%\bibliographystyle{spmpsci}      % mathematics and physical sciences
%\bibliographystyle{spphys}       % APS-like style for physics
%\bibliography{}   % name your BibTeX data base

% Non-BibTeX users please use

\end{document}